\documentclass[11pt]{article}
\usepackage{amsfonts}
\usepackage{amssymb}
\usepackage{color}
\usepackage{amsmath}
=0pt
\def\sqr#1#2{{\vcenter{\vbox{\hrule height.#2pt
              \hbox{\vrule width.#2pt height#1pt \kern#1pt \vrule
width.#2pt}
              \hrule height.#2pt}}}}
\def\dbN{{\mathbb{N}}}

\def\dbR{{\mathbb{R}}}

\def\o{\omega}

\def\3n{\negthinspace \negthinspace \negthinspace }
\def\2n{\negthinspace \negthinspace }
\def\1n{\negthinspace }
\def\ns{\noalign{\smallskip} }

\def\ds{\displaystyle}
\def\D{\Delta}

\def\O{\Omega}
\def\ms{\medskip}

\def\q{\quad}

\def\pa{\partial}

\def\cd{\cdot}

\def\({\Big (}
\def\){\Big )}
\def\[{\Big[}
\def\]{\Big]}

\def\={\buildrel \triangle \over =}

\def\be{\begin{equation}}
\def\bel{\begin{equation}\label}
\def\ee{\end{equation}}
\def\bea{\begin{eqnarray}}
\def\eea{\end{eqnarray}}
\def\bt{\begin{theorem}}
\def\et{\end{theorem}}
\def\bc{\begin{corollary}}
\def\ec{\end{corollary}}
\def\bl{\begin{lemma}}
\def\el{\end{lemma}}
\def\bp{\begin{proposition}}
\def\ep{\end{proposition}}
\def\br{\begin{remark}}
\def\er{\end{remark}}
\def\ba{\begin{array}}
\def\ea{\end{array}}
\def\bd{\begin{definition}}
\def\ed{\end{definition}}

\newtheorem{lemma}{Lemma}[section]
\newtheorem{remark}{Remark}[section]

\newtheorem{theorem}{Theorem}[section]
\newtheorem{corollary}{Corollary}[section]

\newtheorem{definition}{Definition}[section]
\newtheorem{proposition}{Proposition}[section]

\makeatletter
   
   \@addtoreset{equation}{section}
\makeatother

\title{\bf Null control of heat equations with analytic memory kernels}

\author{Qi L\"u\thanks{School of Mathematics,
Sichuan University, Chengdu, 610064, China.
(\tt{lu@scu.edu.cn}).},\q Xu Zhang\thanks{School
of Mathematics, Sichuan University, Chengdu
610064, China. ({\small\tt
zhang$\_$xu@scu.edu.cn}).} \q and \q Enrique
Zuazua\thanks{Departamento de Matem\'aticas,
Universidad Aut\'onoma de Madrid, Cantoblanco,
28049 Madrid -  Spain. (\tt{zuazua@uam.es}). }}

\date{}

\begin{document}

\maketitle

\begin{abstract}
We analyze the control properties of heat
equations with memory terms. We recall
previous results  showing that if the moving
support of the control covers the whole domain
where  heat diffuses, the system is null
controllable when the memory kernel is
polynomial. We formulate the problem of
extending this result to the case of some more
general memory kernels, in particular analytic
ones.

\end{abstract}


\medskip



\section{Formulation of the problem}


Addressed originally in the context of linear
finite-dimensional systems  in \cite{Kalman},
the problem of controllability has been studied
for a broad class of systems including  infinite
dimensional, nonlinear  and stochastic systems.

Many relevant  physical and chemical processes
are effected not only by its current state but
also by its history and the models describing these
processes  involve memory terms.  Some typical
examples are viscoelasticity, non-Fickian
diffusion and thermal processes with memory. We
refer to \cite{AFG, Pruss} and
the rich references therein for more details
about the background and mathematical
formulation of these systems, respectively.

The
model for the control of thermodynamics with memory
and non-Fickian diffusion can be written as follows:
\begin{equation}\label{system4}
\left\{
\begin{array}{ll}\ds
y_{t} -  \int_0^tM(t-s)\D y(s)ds - b\D y =
\chi_{\o(t)}u  &\mbox{ in } (0,T)\times\O,\\[3mm]
\ns\ds y=0 &\mbox{ on } (0,T)\times\pa\O,\\[3mm]
\ns\ds y(0)=y_0  &\mbox{ in }\O.
\end{array}
\right.
\end{equation}
Here $\O\subset\dbR^d$ (for $d\in\dbN$) is a bounded
domain with the $C^2$ boundary $\pa\O$,
$\o(t)\subset\O$ is a nonempty open proper subset for
all $t\in [0,T]$, which stands for the support
of the control that may vary with time, $b\ge
0$, $M(\cd)\in C^2[0,T]$ and $y_0\in L^2(\O)$.

System \eqref{system4} is said to be memory-type null
controllable if for any $y_0\in L^2(\O)$, there
exists a control $u(\cdot)\in L^2(0,T;L^2(\O))$
such that the corresponding solution $y(\cdot)$
satisfies
\begin{equation}\label{eq1}
y(T)=0 \q\hbox{and}\q \int_0^TM(T-s)\D y(s)ds=0.
\end{equation}

The main open problem we formulate is as follows:

\ms

{\bf Problem (P)} {\it  Prove, under suitable
assumptions on the memory kernel $M$ and on the
support $\o(\cdot)$ of the control,
 that the system \eqref{system4} is memory-type
null controllable. }

\ms

One can show that, if $\o(t)$ is
independent of $t$, then the only chance for
obtaining the memory-type null controllability
of \eqref{system4} is to choose $\o=\O$.


\section{Existing results and other related questions}\label{sec_intro}


The above problem can be considered in some
specific situations that lead to some particular
open problems which may deserve separate
analysis.

\begin{itemize}
\item We first consider the
case that $M(\cd)\equiv 1$ and $b=1$. Let
$z(\cd) \= y(\cd)+\int_0^\cd y(s)ds$. Then the
system \eqref{system4} is equivalent to the following one:
\begin{equation}\label{system1}
\left\{
\begin{array}{ll}\ds
z_t-\Delta z -y  = \chi_{\omega(t)}u &  \mbox{ in } (0,T)\times\O,  \\[3mm]
\ns\ds y_t-z_t + y  = 0 &  \mbox{ in }
(0,T)\times\O,
\\[3mm]
z = 0 & \mbox{ on } (0,T)\times\pa\O, \\[3mm]
z(0) = y_0,\;\, y(0)=y_0 & \mbox{ in } \Omega.
\end{array}
\right.
\end{equation}
The null
controllability of \eqref{system1}, ensuring that both $y$ and $z$ can be driven to zero at time $t=T$, is equivalent to the memory-type null controllability of
\eqref{system4}.

For the simpler model in which $z_t$ does not enter into the
second equation of \eqref{system1} the
null controllability of \eqref{system1} was
already established in \cite{CZZ}. However, it
is unclear  whether the technique in
\cite{CZZ}, which is explained below, can be
applied to solve the null controllability of the full
problem of \eqref{system1}.

\vspace{0.2cm}

\item Another simplified version of
\eqref{system4} is the following system:
\begin{equation}\label{heatcc-moving}
\left\{
\begin{array}{ll}
\displaystyle y_t-\Delta y + \int_0^t M(t-s) y(s)ds   = \chi_{\omega(t)}u &  \mbox{ in } (0,T)\times\O,  \\[3mm]
y = 0 & \mbox{ on } (0,T)\times\pa\O, \\[3mm]
y(0) = y_0 & \mbox{ in } \Omega.
\end{array}
\right.
\end{equation}
Here, we replaced the $\int_0^t M(t-s) \D
y(s)ds$ in \eqref{system4} by $\int_0^t M(t-s)
y(s)ds$.
In  \cite{CRZ}, roughly speaking, it is shown that when the
memory kernel $M$ is non-trivial, this system
fails to be null controllable if $\omega$ is a
strict subset of the domain $\Omega$ and
independent of $t$. It is also proved that this
system is null controllable if the support of
the controller $\omega (t)$ is allowed to move, and $\omega(t)$  covers the whole
domain $\Omega$ as $t$ evolves from $t=0$ to
$t=T$ and satisfies some extra geometric
restrictions.

The proof of the above mentioned controllability result is based on viewing
\eqref{heatcc-moving} as the coupling of a heat
equation with an Ordinary Differential Equation
(ODE for short) for the memory term:
$$
z_1(t)\= \int_0^t M(t-s) y(s)ds.
$$
The controlled memory heat equation \eqref{heatcc-moving} can then be written as
\begin{equation}\label{heat-moving}
\left\{
\begin{array}{ll}
\displaystyle y_t-\Delta y + z_1  = u\chi_{\omega(t)}(x) &  \mbox{ in }(0,T)\times\O,  \\[3mm]
\displaystyle z_{1,t}=M(0) + \int_0^t M_t(t-s) y(s)ds &  \mbox{ in }(0,T)\times\O,  \\[3mm]
y = z_1 = 0 & \mbox{ on }(0,T)\times\pa\O, \\[3mm]
y(0) = y_0 , \, z_1(0)=0& \mbox{ in } \Omega.
\end{array}
\right.
\end{equation}
In this way the system \eqref{heatcc-moving} is
reduced to the coupling of a heat equation with
an ODE, the control being applied on the heat
component. But the ODE in \eqref{heat-moving} still involves a memory
term, which is a function of $y$. To cope with
this new memory term we can  introduce a second
auxiliary variable
$$
z_2(t)\= \int_0^t M_t(t-s) y(s)ds.
$$

Iterating this procedure we see that, if $M$ is a polynomial kernel, the controlled memory heat equation can be reduced to a system coupling a heat equation with a finite number of ODEs.

More generally, when the memory kernel $M$ is of the form $
M(t)=e^{at}\sum_{i=0}^K a_i t^i $ for some $K\in\dbN$ and $a,a_i\in\dbR$, the system can be reduced to the coupling of a heat equation with a finite number of ODEs, and the results in
\cite{CRZ} guarantee the null-controllability of
the memory heat equation above.

The proof consists in considering the dual observability inequalities for the adjoint system and deriving them by Carleman inequalities. For this to be done the
condition that $\omega(t)$ covers the whole
domain $\Omega$ in the time interval $[0, T]$ is needed,
together with  the added technical condition that $\Omega$ is split into two disjoint
connected subdomains for all $t$ by $\omega(t)$. This moving geometric condition is necessary to establish the Carleman inequalities for the ODE components of the system.

\end{itemize}

When the kernel $M$ is analytic the same procedure described in the second example leads to a system coupling a heat equation with an infinite number of ODEs.
Whether these techniques, based on Carleman inequalities,  can be  adapted to the
case of analytic kernels is an {\it OPEN}
problem because of the need of dealing with the
superposition of an infinite countable number of
ODEs.

Similar results have been  developed in
\cite{LZZ} for wave equations involving memory
kernels, under suitable moving conditions on the
support of the control $\o(t)$. But, again, the
results in  \cite{LZZ} are limited to  kernels
satisfying suitable conditions.
Whether wave equations with more general smooth
kernels are memory-type controllable is an interesting {\it
OPEN} problem.

\paragraph{Acknowledgments.}

Qi L\"u was supported by the NSF of China under
grant 11471231, and Grant MTM2011-29306-C02-00
of the MICINN, Spain. Xu Zhang was supported by
the NSF of China under grants 11221101 and 11231007, the
PCSIRT under grant IRT1273, the Chang Jiang
Scholars Program from the Chinese Education
Ministry. Enrique Zuazua was partially supported
by the Advanced Grant NUMERIWAVES/FP7-246775 of
the European Research Council Executive Agency,
 the FA9550-15-1-0027 of AFOSR, the
MTM2011-29306 and MTM2014-52347
Grants of the MINECO, and a Humboldt Award at
the University of Erlangen-N\"urnberg. This work was done while the third author was visiting Sichuan University, Chengdu, China.

\end{document}